\documentclass[a4paper]{amsart}

\usepackage[mathscr]{eucal}
\usepackage{amsmath,amsfonts,amssymb}
\usepackage{amsthm}

\newcommand{\mathsym}[1]{{}}

\theoremstyle{plain}
\newtheorem{theorem0}{Theorem}

\newtheorem{theorem}{Theorem}[section]

\newtheorem{lemma}[theorem]{Lemma}

\theoremstyle{definition}
\newtheorem{definition}{Definition}[section]
\newtheorem{remark}{Remark}[section]

\begin{document}
\title[Spherical designs from the norm-$3$ shell of integral lattices]{Spherical designs from norm-$3$ shell of integral lattices}
\author{Junichi Shigezumi}

\maketitle \vspace{-0.1in}
\begin{center}
Graduate School of Mathematics Kyushu University\\
Hakozaki 6-10-1 Higashi-ku, Fukuoka, 812-8581 Japan\\
{\it E-mail address} : j.shigezumi@math.kyushu-u.ac.jp \vspace{-0.05in}
\end{center} \quad

\begin{quote}
{\small\bfseries Abstract.}
A set of vectors all of which have a constant (non-zero) norm value in an Euclidean lattice is called a shell of the lattice. Venkov classified strongly perfect lattices of minimum $3$ (R\'{e}seaux et ``designs'' sph\'{e}rique, 2001), whose minimal shell is a spherical $5$-design. This note considers the classification of integral lattices whose shells of norm $3$ are $5$-designs.\\  \vspace{-0.15in}

\vspace{-0.15in}

 \vspace{0.15in}
\end{quote}

\section*{Introduction}

Let $L$ be an Euclidean lattice, which is a discrete vector space over $\mathbb{Z}$. The squared norm of a vector of the lattice is called the {\it norm of the vector}. Then, the set $s_m (L)$ of vectors of the lattice $L$ which take the same value $m$ for their norm is called the {\it shell of the lattice}, i.e. $s_m (L) := \{ x \in L \: ; \: (x, x) = m \}$. Moreover, the shell of {\it minimum} $\min_{x \in L \setminus \{ 0 \}} (x, x)$ of the lattice $L$ is called the {\it minimal shell}, which is denoted by $S(L)$.

\begin{definition}[Spherical design \cite{DGS}]\label{def-design}
Let $X$ be a non-empty finite set on the Euclidean sphere $S^{d-1}$, and let $t$ be a positive integer. $X$ is called a spherical $t$-design if
\begin{equation}
\frac{1}{| S^{d-1} |} \int_{S^{d-1}} f(\xi) \, d \xi \ = \ \frac{1}{| X |} \hspace{0.05in} \sum_{\xi \in X} \hspace{0.05in} f(\xi) \label{eq-design}
\end{equation}
for every polynomial $f(x) = f(x_1, \ldots, x_d)$ of degree at most $t$.
\end{definition}

For every nonempty shell $s_m (L)$ of a lattice $L$, a normalization $X = \frac{1}{\sqrt{m}} s_m (L)$ is considered, where $X$ is a finite set on an Euclidean sphere. A lattice, whose minimal shell is a spherical $4$-design (i.e. a $5$-design), is said to be {\it strongly perfect}.

B. B. Venkov proved the following theorem:

\begin{theorem0}[Venkov \cite{V}, Theorem 7.4]\label{th-ven}
The strongly perfect lattices that are integral and of minimum $3$ are $O_1$, $O_7$, $O_{16}$, $O_{22}$, and $O_{23}$. Furthermore, the minimal shell is a spherical $7$-design only for the case of the lattice $O_{23}$.
\end{theorem0}

Now, as an expansion of the above theorem, we prove the following theorem:

\begin{theorem0}\label{th-main}
Let $L$ be an integral lattice. If its shell of norm $3$ is a spherical $5$-design, then $L$ is isometric to one of the following nine lattices:
\begin{enumerate}
\item $\mathbb{Z}^7$, whose minimum is equal to one.

\item $\Lambda_{16, 2, 1}$, $\Lambda_{16, 2, 2}$ and $\Lambda_{16, 2, 3}$, whose minima are equal to two.

\item $O_1$, $O_7$, $O_{16}$, $O_{22}$, and $O_{23}$, whose minima are equal to three.
\end{enumerate}
\end{theorem0}

The definitions of the lattices in the above theorems are given in the next section. The remaining sections are devoted to the proof of Theorem \ref{th-main}.\\

\section{Definitions of the lattices in Theorem \ref{th-main}}\label{sec-def}

\subsection{Lattices of minimum $3$}

\begin{lemma}[Venkov \cite{V}, Lemma 7.1]
Let $L$ be an even integral lattice of dimension $n \geq 2$ and of minimum $4$, and let $e$ be a minimal vector of $L$. Denote by $p$ the orthogonal projection on the hyperplane $H = e^{\perp}$, put ${L_e}' = \{ x \in L \ | \ (e, x) \equiv 0 \pmod{2} \}$, and let $L_e = p({L_e}')$. Suppose that one of the following two assumptions holds:
\def\labelenumi{(\arabic{enumi})}
\begin{enumerate}
\item There is $x \in L$ such that $(e, x) \equiv 1 \pmod{2}$;

\item We have $(y, e) \equiv 0 \pmod{2}$ for all $y \in L$, and $L$ contains a vector $x$ such that $(e, x) \equiv 2 \pmod{4}$.
\end{enumerate}
Then, $L_e$ is a odd integral lattice of minimum at least $3$, and we have $\det(L_e) = \det(L)$ under assumption $(1)$ and $\det(L_e) = \frac{1}{4} \det(L)$ under assumption $(2)$.
\end{lemma}

We denote by $\Lambda_n$ the {\it laminated lattices} for $2 \leq n \leq 24$ (See Conway-Sloane \cite{CS}, Ch. 6). Note that $\Lambda_n$ is isometric to $\sqrt{2} \mathbb{E}_n$ for $n = 6$, $7$, $8$, that $\Lambda_{16}$ is isometric to the {\it Barnes-Wall lattice} $BW_{16}$, and that $\Lambda_{24}$ is the {\it Leech lattice}. Then, we set $O_1 = \sqrt{3} \mathbb{Z}$. We denote by $O_7$ (resp. $O_{23}$) the projected $\L_e$ associated with the laminated lattice $\Lambda_8$ (resp. $\Lambda_{24}$). Finally, we denote by $O_{22}$ (resp. $O_{16}$) the orthogonal of $O_1$ (resp. $O_7$) in $O_{23}$.

We have $\det (O_1) = \det (O_{22}) = 3$, $\det (O_7) = \det (O_{16}) = 64$, and $\det (O_{23}) = 1$; thus $O_{23}$ is unimodular.

The theta series of each lattice have the following form:
\begin{align*}
\Theta_{O_1} &= 1 + 2 \, q^3 + 2 \, q^{12} + 2 \, q^{27} + 2 \, q^{48} + \cdots \\
\Theta_{O_7} &= 1 + 56 \, q^3 + 126 \, q^4 + 576 \, q^7 + 756 \, q^8 + 1512 \, q^{11} + 2072 \, q^{12}\\
 &+ 4032 \, q^{15} + 4158 \, q^{16} + 5544 \, q^{19} + 7560 \, q^{20} + 12096 \, q^{23} + 11592 \, q^{24} + \cdots \\
\Theta_{O_{16}} &= 1 + 512 \, q^3 + 4320 \, q^4 + 18432 \, q^5 + 61440 \, q^6 + 193536 \, q^7\\
 &+ 522720 \, q^8 + 1126400 \, q^9 + 2211840 \, q^{10} + 4584960 \, q^{11} + 8960640 \, q^{12} + \cdots \\
\Theta_{O_{22}} &= 1 + 2816 \, q^3 + 49896 \, q^4 + 456192 \, q^5 + 2821632 \, q^6 + 13229568 \, q^7\\
 &+ 50332590 \, q^8 + 163175936 \, q^9 + 467596800 \, q^{10} + 1214196480 \, q^{11}\\
 &+ 2900976144 \, q^{12} + \cdots \\
\Theta_{O_{23}} &= 1 + 4600 \, q^3 + 93150 \, q^4 + 953856 \, q^5 + 6476800 \, q^6 + 32788800 \, q^7\\
 &+ 133204500 \, q^8 + 458086400 \, q^9 + 1384998912 \, q^{10} + 3771829800 \, q^{11}\\
 &+ 9403968600 \, q^{12} + \cdots
\end{align*}
Since $O_{23}$ is unimodular, we also have that $\Theta_{O_{23}} = \theta_3^{23} - 46 \, \theta_3^{15} \Delta_8$, where $\Delta_8 = \frac{1}{16} \theta_2^4 \theta_4^4$ and $\theta_i$ for $i = 2, 3, 4$ are known as {\it Jacobi's theta functions}. (See \cite{CS}, \cite{N})

Let $X$ be a nonempty finite set on the Euclidean sphere $S^{d-1}$ ($\subset \mathbb{R}^{d}$). We denote the distance set of $X$ by $A(X) := \{ (x, y) \: ; \: x, y \in X, x \ne y \}$; then we call $X$ an $s$-distance set if $| A(X) | = s$. Now, $X$ is said to be a $(d, n, s, t)$-configuration if $X \subset S^{d-1}$ is of order $n (:= | X |)$, a $s$-distance set, and a spherical $t$-design. The following table contains the $(d, n, s, t)$-configuration of each shell of norm $m$ of the lattice:
\begin{center}
\begin{tabular}{r}
\quad\\
\hline
$m$\\
\hline
$3$\\
$4$\\
$5$\\
$6$\\
$7$\\
$8$\\
$9$\\
$10$\\
$11$\\
$12$\\
\hline
\end{tabular}
\quad
\begin{tabular}{cccc}
\multicolumn{4}{c}{$O_7$}\\
\hline
$d$ & $n$ & $s$ & $t$\\
\hline
$7$ & $56$ & $3$ & $5$\\
$7$ & $126$ & $4$ & $5$\\
\ & \ & \ & \\
\ & \ & \ & \\
$7$ & $576$ & $7$ & $5$\\
$7$ & $756$ & $8$ & $53$\\
\ & \ & \ & \\
\ & \ & \ & \\
$7$ & $1512$ & $11$ & $5$\\
$7$ & $2072$ & $12$ & $5$\\
\hline
\end{tabular}
\quad
\begin{tabular}{cccc}
\multicolumn{4}{c}{$O_{16}$}\\
\hline
$d$ & $n$ & $s$ & $t$\\
\hline
$16$ & $512$ & $4$ & $5$\\
$16$ & $4320$ & $6$ & $7$\\
$16$ & $18432$ & $8$ & $5$\\
$16$ & $61440$ & $10$ & $7$\\
$16$ & $193536$ & $12$ & $5$\\
$16$ & $522720$ & $14$ & $7$\\
\hline
\quad\\ \quad\\ \quad\\ \quad\\
\end{tabular}

\quad\\

\begin{tabular}{r}
\quad\\
\hline
$m$\\
\hline
$3$\\
$4$\\
$5$\\
\hline
\end{tabular}
\quad
\begin{tabular}{cccc}
\multicolumn{4}{c}{$O_{22}$}\\
\hline
$d$ & $n$ & $s$ & $t$\\
\hline
$22$ & $2816$ & $4$ & $5$\\
$22$ & $49896$ & $6$ & $5$\\
$22$ & $456192$ & $8$ & $5$\\
\hline
\end{tabular}
\quad
\begin{tabular}{cccc}
\multicolumn{4}{c}{$O_{23}$}\\
\hline
$d$ & $n$ & $s$ & $t$\\
\hline
$23$ & $4600$ & $4$ & $7$\\
$23$ & $93150$ & $6$ & $7$\\
$23$ & $953856$ & $8$ & $7$\\
\hline
\end{tabular}
\end{center}\quad

\subsection{Lattices of minimum $2$}

Let $\varepsilon_1, \ldots, \varepsilon_{16}$ be an orthonormal basis of $\mathbb{R}^{16}$. We denote some vectors
\begin{gather*}
f_1 := \frac{\varepsilon_1 + \cdots + \varepsilon_8}{2} + \varepsilon_9, \quad
f_2 := \varepsilon_1 + \frac{\varepsilon_9 + \cdots + \varepsilon_{16}}{2},\\
f_3 := \varepsilon_1 + \varepsilon_5 + \varepsilon_9 + \varepsilon_{13},\\
f_4 := \varepsilon_1 + \varepsilon_3 + \varepsilon_5 + \varepsilon_7, \quad
f_5 := \varepsilon_1 + \varepsilon_3 + \varepsilon_9 + \varepsilon_{11}, \quad
f_6 := \varepsilon_1 + \varepsilon_3 + \varepsilon_{13} + \varepsilon_{15}.
\end{gather*}

Now, we define the following three lattices
\begin{align}
\Lambda_{16,2,1} &:= \langle (A_1)^{16}, f_1, f_2, f_3, f_4, f_5, f_6 \rangle, \\
\Lambda_{16,2,2} &:= \langle (D_4)^4, f_1, f_2, f_3 \rangle, \\
\Lambda_{16,2,3} &:= \langle (D_8)^2, f_1, f_2 \rangle,
\end{align}
where we put root systems in the above definitions as $(A_1)^{16} := \{ \pm (\varepsilon_{2 i - 1} \pm \varepsilon_{2 i}) \: ; \: 1 \leqslant i \leqslant 8 \}$, $(D_4)^4 := \{ \pm (\varepsilon_i \pm \varepsilon_j) \: ; \: 1 \leqslant i < j \leqslant 4, 5 \leqslant i < j \leqslant 8, 9 \leqslant i < j \leqslant 12, \: \text{or} \; 13 \leqslant i < j \leqslant 16 \}$, and $(D_8)^2 := \{ \pm (\varepsilon_i \pm \varepsilon_j) \: ; \: 1 \leqslant i < j \leqslant 8 \; \text{or} \; 9 \leqslant i < j \leqslant 16 \}$. Then, we have $(A_1)^{16} \subset (D_4)^4 \subset (D_8)^2$ and $\Lambda_{16,2,1} \subset \Lambda_{16,2,2} \subset \Lambda_{16,2,3}$. Furthermore, we have
\begin{equation}
\Lambda_{16,2,3}  = \Lambda_{16,2,2} \cup (\varepsilon_1 + \varepsilon_5 + \Lambda_{16,2,2})
\quad \text{and} \quad
\Lambda_{16,2,2}  = \Lambda_{16,2,1} \cup (\varepsilon_1 + \varepsilon_3 + \Lambda_{16,2,1}).
\end{equation}

\begin{remark}
We denote some other vectors
\begin{align*}
f_7 &:= \varepsilon_1 + \varepsilon_2 + \varepsilon_3 + \varepsilon_4, \quad
f_8 := \varepsilon_1 + \varepsilon_2 + \varepsilon_5 + \varepsilon_6, \quad
f_9 := \varepsilon_1 + \varepsilon_2 + \varepsilon_7 + \varepsilon_8, \\
f_{10} &:= \varepsilon_1 + \varepsilon_2 + \varepsilon_9 + \varepsilon_{10}, \quad
f_{11} := \varepsilon_1 + \varepsilon_2 + \varepsilon_{11} + \varepsilon_{12}, \quad
f_{12} := \varepsilon_1 + \varepsilon_2 + \varepsilon_{13} + \varepsilon_{14}, \\
f_{13} &:= \varepsilon_1 + \varepsilon_2 + \varepsilon_{15} + \varepsilon_{16}.
\end{align*}
Then, we can write
\begin{equation*}
O_{16} = \langle (\sqrt{2} A_1)^{16}, f_1, f_2, f_3, f_4, f_5, f_6, f_7, f_8, f_9, f_{10}, f_{11}, f_{12}, f_{13} \rangle,
\end{equation*}
and we have $(\sqrt{2} A_1)^{16} := \{ \pm 2 \varepsilon_i \: ; \: 1 \leqslant i \leqslant 16 \} \subset (A_1)^{16}$ and $O_{16} \subset \Lambda_{16,2,1}$. Furthermore, we have
\begin{equation}
\Lambda_{16,2,1}  = O_{16} \cup (\varepsilon_1 + \varepsilon_2 + O_{16}).
\end{equation}
\end{remark}\quad

We have $\det (\Lambda_{16,2,1}) = 16$, $\det (\Lambda_{16,2,2}) = 4$, and $\det (\Lambda_{16,2,3}) = 1$; thus $\Lambda_{16,2,3}$ is unimodular.

We obtain the theta series of the lattices by numerical calculation as the following form:
\begin{align*}
\Theta_{\Lambda_{16,2,1}} &= 1 + 32 \, q^2 + 1024 \, q^3 + 8160 \, q^4 + 36864 \, q^5 + 127360 \, q^6 + 387072 \, q^7\\
 &+ 1016288 \, q^8 + 2252800 \, q^9 + 4564416 \, q^{10} + 9169920 \, q^{11} + 17395328 \, q^{12} + \cdots \\
\Theta_{\Lambda_{16,2,2}} &= 1 + 96 \, q^2 + 2048 \, q^3 + 15840 \, q^4 + 73728 \, q^5 + 259200 \, q^6 + 774144 \, q^7\\
 &+ 2003424 \, q^8 + 4505600 \, q^9 + 9269568 \, q^{10} + 18339840 \, q^{11} + 34264704 \, q^{12} +  \cdots \\
\Theta_{\Lambda_{16,2,3}} &= 1 + 224 \, q^2 + 4096 \, q^3 + 31200 \, q^4 + 147456 \, q^5 + 522880 \, q^6 + 1548288 \, q^7\\
 &+ 3977696 \, q^8 + 9011200 \, q^9 + 18679872 \, q^{10} + 36679680 \, q^{11} + 68003456 \, q^{12} +  \cdots 
\end{align*}
Since $\Lambda_{16,2,3}$ is unimodular, we also have that $\Theta_{\Lambda_{16,2,3}} = \theta_3^{16} - 32 \theta_3^8 \Delta_8$.

The following table is the $(d, n, s, t)$-configuration of each shell of norm $m$ of the lattice:
\begin{center}
\begin{tabular}{r}
\quad\\
\hline
$m$\\
\hline
$2$\\
$3$\\
$4$\\
$5$\\
$6$\\
$7$\\
$8$\\
$9$\\
\hline
\end{tabular}
\quad
\begin{tabular}{cccc}
\multicolumn{4}{c}{$\Lambda_{16,2,1}$}\\
\hline
$d$ & $n$ & $s$ & $t$\\
\hline
$16$ & $32$ & $2$ & $3$\\
$16$ & $1024$ & $6$ & $5$\\
$16$ & $8160$ & $8$ & $3$\\
$16$ & $36864$ & $10$ & $5$\\
$16$ & $127360$ & $12$ & $3$\\
$16$ & $387072$ & $14$ & $5$\\
$16$ & $1016288$ & $16$ & $3$\\
$16$ & $2252800$ & $18$ & $5$\\
\hline
\end{tabular}
\quad
\begin{tabular}{cccc}
\multicolumn{4}{c}{$\Lambda_{16,2,2}$}\\
\hline
$d$ & $n$ & $s$ & $t$\\
\hline
$16$ & $96$ & $4$ & $3$\\
$16$ & $2048$ & $6$ & $5$\\
$16$ & $15840$ & $8$ & $3$\\
$16$ & $73728$ & $10$ & $5$\\
$16$ & $259200$ & $12$ & $3$\\
$16$ & $774144$ & $14$ & $5$\\
$16$ & $2003424$ & $16$ & $3$\\
\hline
\quad\\
\end{tabular}
\quad
\begin{tabular}{cccc}
\multicolumn{4}{c}{$\Lambda_{16,2,3}$}\\
\hline
$d$ & $n$ & $s$ & $t$\\
\hline
$16$ & $224$ & $4$ & $3$\\
$16$ & $4096$ & $6$ & $5$\\
$16$ & $31200$ & $8$ & $3$\\
$16$ & $147456$ & $10$ & $5$\\
$16$ & $522880$ & $12$ & $3$\\
$16$ & $1548288$ & $14$ & $5$\\
\hline
\quad\\ \quad\\
\end{tabular}
\end{center}

\subsection{Lattice of minimum $1$}

We have $\det (\mathbb{Z}^7) = 1$, thus $\mathbb{Z}^7$ is unimodular.

The theta series of the lattices have the following form:
\begin{align*}
\Theta_{\mathbb{Z}^7} &= 1 + 14 \, q + 84 \, q^2 + 280 \, q^3 + 574 \, q^4 + 840 \, q^5 + 1288 \, q^6 + 2368 \, q^7\\
 &+ 3444 \, q^8 + 3542 \, q^9 + 4424 \, q^{10} + 7560 \, q^{11} + 9240 \, q^{12} + \cdots
\end{align*}
Since $\mathbb{Z}^7$ is unimodular, we also have $\Theta_{\mathbb{Z}^7} = \theta_3^7$.

For the spherical design from each shell of $\mathbb{Z}^7$, the following facts are already known:
\begin{theorem}[Pache \cite{P}, parts of Theorem 25 and Proposition 26]\quad
\begin{enumerate}
\item For $n \geqslant 2$, all the nonempty shells of $\mathbb{Z}^n$ are spherical $3$-designs.

\item The following shells are spherical $5$-designs:
\begin{align*}
&s_m (\mathbb{Z}^4) & &m = 2 a, \quad a \geqslant 1.\\
&s_m (\mathbb{Z}^7) & &m = 4^a (8 b + 3), \quad a, b \geqslant 0.
\end{align*}

\item For $n \geqslant 2$ and $1 \leqslant m \leqslant 1200$, the nonempty shells of norm $m$ of \ $\mathbb{Z}^n$ are not spherical $5$-designs, except for the above cases.
\end{enumerate}
\end{theorem}

\begin{remark}
$\mathbb{Z}^7$ and $O_7$ have $8 \, \binom{7}{3} = 280$ and $8 \cdot 7 = 56$ vectors of norm $3$, respectively. Then, as a natural question, can we write $s_3(\mathbb{Z}^7)$ as a disjoint union of configurations isometric to $s_3(O_7)$? The answer is no.\\
There are $30$ subsets of $s_3(\mathbb{Z}^7)$ which are isometric to $s_3(O_7)$. However, any $3$ such subsets are not disjoint. Here, we can choose $2$ disjoint subsets, for example, $(\pm 1, \pm 1, 0, \pm 1, 0, 0, 0)^C$ and $(\pm 1, 0, \pm 1, \pm 1, 0, 0, 0)^C$, where ``$\pm$'' indicates that we take all possible sign changes, and $C$ indicates that we take any cyclic shifts.
\end{remark}

\begin{remark}
Note that all the lattices in this section are $3$-lattices, which are generated by some vectors of norm $3$.
\end{remark}\quad

Let $L$ be an integral lattice of dimension $n$, whose shell of norm $3$ is a spherical $5$-design. By $X := s_3 (L)$ the shell of norm $3$ is denoted. The argument of Theorem \ref{th-main} $(3)$ is just equivalent to Theorem \ref{th-ven}. Thus, we may suppose $\min (L)$ equal to $1$ or $2$.\\

\section{On spherical designs}

\begin{theorem}[Venkov \cite{V}, Theorem 3.1]
Let $X  \in S^{n-1}$ be a finite set, and $t$ be a positive integer. By $e$ $($resp. $o$ $)$ the greatest even $($resp. odd $)$ integer which is at most $t$ is denoted. Then, $X$ is a spherical $t$-design if and only if there is a constant $c_e$ such that, for every $\alpha \in \mathbb{R}^n$, we have the two equations
\begin{equation}
\sum_{x \in X} (x, \alpha)^e = c_e (\alpha, \alpha)^{e/2} \quad \text{and} \quad \sum_{x \in X} (x, \alpha)^o= 0.
\end{equation}
\end{theorem}

If the above two equations hold, by repetition of the Laplacian $\Delta_{\alpha}$, we always have the formulae
\begin{equation*}
\sum_{x \in X} (x, \alpha)^k = c_k (\alpha, \alpha)^{k/2} \quad \text{and} \quad \sum_{x \in X} (x, \alpha)^l = 0
\end{equation*}
for any even $k \leq e$ and any odd $l \leq o$, where the notation $\Delta_y$ refers to derivation with respect to the variable $y$. We also have
\begin{equation*}
c_k = \frac{1 \cdot 3 \cdot 5 \cdots (k-1)}{n (n+2) \cdots (n + k - 2)} |X|.
\end{equation*}

In this paper, the finite sets on spheres from shells of Euclidean lattices are considered. Then, every set is antipodal, thus the second equation always holds. Also, $5$-designs from vectors of norm $3$ are also considered here. Thus it is necessary and sufficient to consider the following two equations:
\begin{align}
\sum_{x \in X} (x, \alpha)^2 &= \frac{3 |X|}{n} (\alpha, \alpha), \label{eq-sum2}\\
\sum_{x \in X} (x, \alpha)^4 &= \frac{27 |X|}{n (n+2)} (\alpha, \alpha)^2. \label{eq-sum4}
\end{align}\quad

Again, let $X$ be the shell of norm $3$ of the lattice $L$. For any vectors $x_0 \in X$, we denote $n_i := | \{ x \in X \: ; \: (x_0, x) = i \} |$ for $i = 0, 1, 2$. By the above equations, taking $\alpha = x_0$, we have
\begin{equation}
n_0 = \frac{4 n^2 -37 n +153}{4 n (n + 2)} |X| - 20, \quad
n_1 = \frac{3 (4 n - 19)}{2 n (n + 2)} |X| + 15, \quad
n_2 = \frac{3 (25 - n)}{8 n (n + 2)} |X| - 6. \label{eq-ni0}
\end{equation}
Note that these results do not depend on the choice of $x_0$.\\

\section{Minimum of lattices}

Let $t \in S(L)$ be a minimal vector of the lattice $L$. Since $(x \pm t, x \pm t) - (t, t) = (x, x) \pm 2 (x, t) \geq 0$ for any $x \in X$, we have $|(x, t)| \leqslant \frac{1}{2} (x, x) = \frac{3}{2}$ (cf. \cite{V}, Lemma 6.10). Thus, we have
\begin{equation}
(x, t) \in \{ 0, \pm 1 \}. \label{bd-inpro}
\end{equation}
Then, we denote $p_i := | \{ x \in X \: ; \: (x, t) = i \} |$ for $i = 0, 1$. By the equalities (\ref{eq-sum2}) and (\ref{eq-sum4}), taking $\alpha = t$, we have
\begin{equation}
(t, t) = \frac{n + 2}{9}, \quad
p_0 = \frac{2(n - 1)}{3 n} |X|, \quad
p_1 = \frac{n + 2}{6 n} |X|.\label{eq-pi0}
\end{equation}
(cf. \cite{V}, Lemma 7.11)\\

If $\min (L) = 1$, by the first equation, the dimension of the lattice $L$ is equal to $7$. By \cite{S}, there are only two $3$-lattices whose shells of norm $3$ are spherical $5$-designs, which are $O_7$ and $\mathbb{Z}^7$. Here, a $3$-lattice is an integral lattice which is generated by vectors of norm $3$. If we consider the lattice $L'$ which is generated by $s_3 (L)$, then $L'$ is a $3$-lattice such that $s_3 (L') = s_3 (L)$. Thus, $|X|$ is equal to $56$ or $280$.

Let $|X| = 56$; then $n_0 = n_2 = 0$ by the equality (\ref{eq-ni0}). Since $p_0 > 0$ and $p_1 > 0$, we can take some elements $x1, x_2 \in X$ such that $(x_1, t) = 0$, $(x_2, t) = 1$, and $(x_1, x_2) = \pm 1$. Then, $x_1 \mp (x_2 - t)$ is a vector of norm $3$, i.e. $x_1 \mp (x_2 - t) \in X$, and we have $(x_1 \mp (x_2 - t), x_1) = 2$. This contradicts $n_2 = 0$.

Let $|X| = 280$; then we have $L \supset L' \simeq \mathbb{Z}^7$. Since $L$ is integral, we have $L \simeq \mathbb{Z}^7$.

In conclusion, if $\min (L) = 1$, then we have $L \simeq \mathbb{Z}^7$.\\

Now, the remaining case is when $\min (L) = 2$.

If $\min (L) = 2$, by the equations (\ref{eq-ni0}) and (\ref{eq-pi0}), we have the following equations:
\begin{gather}
n = 16, \quad
p_0 = \frac{5}{8} |X|, \quad
p_1 = \frac{3}{16} |X|.\label{eq-pi}\\
n_0 = \frac{65}{128} |X| - 20, \quad
n_1 = \frac{15}{64} |X| + 15, \quad
n_2 = \frac{3}{256} |X| - 6. \label{eq-ni}
\end{gather}
Since $n_i$ and $p_i$ are nonnegative integers,
\begin{equation}\label{eq-s3}
256 \, \big| \, |X| \quad \text{and} \quad |X| \geqslant 512.
\end{equation}

Furthermore, for any $x_0 \in X$, if $x \in X$ satisfies $(x_0, x) = 2$, then $x_0 - x \in s_2 (L)$ and $(x_0 - x, x_0) = 1$. On the other hand, if $y \in s_2 (L)$ satisfies $(x_0, y) = 1$, then $x_0 - y \in X$ and $(x_0 - y, x_0) = 2$. Thus, we have
\begin{equation}\label{eq-m1}
|\{ y \in s_2 (L) \: ; \: (x_0, y) = 1 \}| = |\{ x \in X \: ; \: (x_0, x) = 2 \}| = n_2,
\end{equation}
where this number does not depend on the choice of $x_0$. Then, we have
\begin{align*}
n_2 \, |X| &= |\{ y \in s_2 (L) \: ; \: (x_0, y) = 1 \}| \times |X|\\
 &= |\{ x \in X \: ; \: (x, y_0) = 1 \}| \times |s_2 (L)| \, = \, p_1 \, |s_2 (L)|.
\end{align*}
Thus, we have
\begin{equation}\label{eq-s20}
|s_2 (L)| = \frac{1}{16} |X| - 32.
\end{equation}\quad

\section{Intersection numbers}

Let $\alpha, \beta, \gamma \in A(X)$, where $A(X) = \{ (x, y) ; x, y \in X, x \ne y \}$ is a distance set. Then, we choose a pair of vectors $x, y \in X$ such that $(x, y) = \gamma$, and denote
\begin{equation}
P_{\gamma} (\alpha, \beta) := |\{ z \in X \: ; \: (x, z) = \alpha, (z, y) = \beta \}|.
\end{equation}
If this number is uniquely determined for any choice of the pair $x, y$, then it can be called the {\it intersection number}. Since $X$ is antipodal, we have $P_{\gamma} (\alpha, \beta) = P_{\gamma} (\beta, \alpha)$ and $P_{\gamma} (\alpha, \beta) = P_{\gamma} (- \alpha, - \beta) = P_{- \gamma} (\alpha, - \beta)$.

In this section, the intersection numbers  $P_2 (\alpha, \beta)$ are considered. Now, a pair of vectors $x, y \in X$ such that $(x, y) = 2$ is chosen. Then, we have $x - y \in s_2 (L)$. For any $z \in X$, we have $(x - y, z) = (x, z) - (y, z) \in \{ 0, \pm 1\}$ by the relation (\ref{bd-inpro}). Here, we have $P_2 (\alpha, \beta) = 0$ for every $\alpha, \beta$ such that $|\alpha - \beta| > 1$. Furthermore, it is clear that $P_2 (3, 3) = 0$ and $P_2 (2, 3) = 1$. This can be  denoted
\begin{equation*}
a_1 := P_2 (2, 2), \quad a_2 := P_2 (1, 2), \quad a_3 := P_2 (1, 1), \quad a_4 := P_2 (0, 1), \quad a_5 := P_2 (0, 0).
\end{equation*}
We have
\begin{equation*}
n_0 = 2 a_4 + a_5, \quad n_1 = a_2 + a_3 + a_4, \quad n_2 = 1 + a_1 + a_2.
\end{equation*}
By the equations (\ref{eq-sum2}) and (\ref{eq-sum4}), taking $\alpha = x + y$ and $\alpha = x - y$, we obtain four relations for $a_1, \ldots, a_5$. Finally, we have
\begin{gather*}
P_2 (2, 2) = \frac{1}{256} |X| - 4, \quad P_2 (1, 2) = \frac{1}{128} |X| - 3, \quad P_2 (1, 1) = \frac{9}{64} |X| +16,\\
P_2 (0, 1) = \frac{11}{128} |X| + 2, \quad P_2 (0, 0) = \frac{43}{128} |X| - 24.
\end{gather*}

Now, for any element $y_0 \in s_2 (L)$, there exists $x_0 \in X$ such that $(x_0, y_0) = 1$. Then, we have $x_0 - y_0 \in X$ and $(x_0 - y_0, x_0) = 2$. Let $y \in s_2 (L)$ satisfy $(y_0, y) = 1$, then $y_0 - y \in s_2 (L)$ and $(x_0, y) = 0$ or $1$, because $(x_0, y_0 - y) = 1 - (x_0, y) \in \{ 0, \pm 1 \}$ by (\ref{bd-inpro}). If $(x_0, y) = 1$, then $x_0 - y \in X$ and $(x_0, x_0 - y) = (x_0 - y_0, x_0 - y) = 2$. Thus, we have
\begin{equation*}
|\{ y \in s_2 (L) \: ; \: (y_0, y) = 1, (x_0, y) = 1 \}| = |\{ x \in X \: ; \: (x_0, x) = (x_0 - y_0, x) = 2 \}| = P_2 (2, 2).
\end{equation*}
On the other hand, if $(x_0, y) = 0$, then $y_0 - y \in X$ and $(y_0, y_0 - y) = (x_0, y_0 - y) = 1$. Thus, we have
\begin{equation*}
|\{ y \in s_2 (L) \: ; \: (y_0, y) = 1, (x_0, y) = 0 \}| = |\{ y \in s_2 (L) \: ; \: (y_0, y) = 1, (x_0, y) = 1 \}| = P_2 (2, 2).
\end{equation*}
In conclusion, we have
\begin{equation}\label{eq-s21}
|\{ y \in s_2 (L) \: ; \: (y_0, y) = 1 \}| = \frac{1}{128} |X| - 8, \quad |\{ y \in s_2 (L) \: ; \: (y_0, y) = 0 \}| = \frac{1}{64} |X| - 18.
\end{equation}
These relations imply that $s_2 (L)$ is a spherical $3$-design.

In addition, we have the following fact:
\begin{lemma}[see \cite{P}]
Let $L$ be an integral lattice. Then, its shell $s_2 (L)$ of norm $2$ is a root system.
\end{lemma}
Furthermore, irreducible root systems have been classified; they are $A_n$ for $n \geqslant 1$, $D_n$ for $n \geqslant 4$, and $E_n$ for $n = 6, 7, 8$ (see \cite{B}). Orthogonal unions of irreducible root systems which satisfy the conditions (\ref{eq-s3}), (\ref{eq-s20}) and (\ref{eq-s21}) result only in the following nine cases:
\begin{equation}
s_2 (L) \simeq (A_1)^{16}, \, (A_2)^8, \, (A_4)^4, \, (A_8)^2, \, A_{16}, \,
 (D_4)^4, \, (D_8)^2, \, D_{16}, \; \text{or} \; (E_8)^2.
\end{equation}
In the next section, each case is examined. Finally, only the three cases $s_2 (L) = (A_1)^{16}$, $(D_4)^4$, and $(D_8)^2$ will have to be considered, and obtain lattices $\Lambda_{16, 2, 1}$, $\Lambda_{16, 2, 2}$, and $\Lambda_{16, 2, 3}$, respectively.\\

\section{Classification of lattices}

\subsection{The case of $s_2 (L) = (A_1)^{16}$}\label{sec-16-2-1}\quad

We put $(A_1)^{16} = \{ \pm \sqrt{2} \, e_i \: ; \: 1 \leqslant i \leqslant 16 \}$ with an orthonormal basis $\{ e_i \}_{1 \leqslant i \leqslant 16}$ of $\mathbb{R}^{16}$. Consider a vector $x_0 \in X$. We write $x_0 = (a_1, \ldots, a_{16}) \in X$, then we have $(x, \sqrt{2} \, e_i) = \sqrt{2} \, a_i \in \{ 0, \pm 1 \}$, thus we have $a_i \in \{ 0, \pm 1 / \sqrt{2} \}$. In addition, we have $(x, x) = a_1^2 + \cdots + a_{16}^2 = 3$ by definition. Thus, six coordinates are $\pm 1 \sqrt{2}$ and ten coordinates are $0$.

If $a_i = 1 / \sqrt{2}$, then the $i$th coordinate of $x_0 - \sqrt{2} \, e_i$ is equal to $- 1 / \sqrt{2}$. Thus, in the lattice which is an additional group, we take all possible sign changes of nonzero coordinates of $x_0$. We write $\overline{x_0} = \sqrt{2} (|a_1|, \ldots, |a_{16}|)$, then we define equivalence classes $\overline{X} := \{ \overline{x} \: ; \: x \in X \}$, where we regard $\overline{x}$ as a equivalence class. Each class has $2^6$ vectors by sign changes. Since $|X| = 1024$, we have $16$ classes.

Take $x_1 \not\in \overline{x_0}$, and let $l := (\overline{x_0}, \overline{x_1})$. In this case, we have $n_0 = 500$, $n_1 = 255$, and $n_2 = 6$. Firstly, since $(x_1, x_0) \in \{ 0, \pm 1, \pm 2 \}$, we have $l = 0, 2, 4$. We write ${m_i}' := |\{ x \in \overline{x_0} \: ; \: (x, x_0) = i \}|$ and ${m_{l, i}}'' := |\{ x \in \overline{x_1} \: ; \: (x, x_0) = i \}|$ for $i = 0, 1, 2$. Then, we have ${m_0}' = 20$, ${m_1}' = 15$, and ${m_2}' = 6$. Secondly, since ${m_2}' = n_2$, we need ${m_{l, 2}}'' = 0$, thus we have $l = 0, 2$. Finally, we have ${m_{0, 0}}'' = 64$, ${m_{2, 0}}'' = 32$, and $n_0 = 20 + 32 \times 15$ for $16$ equivalent classes, so we have $l = 2$.

Now, if we regard $\overline{X}$ as an incidence matrix of a block design (cf. matrix (\ref{des-inc})) , we can consider a $2$-$(16, 6, 2)$ ($t$-$(v, k, \lambda)$) design. By Gibbons \cite{G}, we have just three equivalence classes for this block design. Furthermore, if we regard the three classes as a basis for a linear code of $\mathbb{F}_2^{16}$, then we obtain three linear codes whose parameters are respectively $[v, k, d] = [16, 6, 6]$, $[16, 7, 4]$, and $[16, 8, 4]$. Moreover, when $[v, k, d] = [16, 7, 4]$ and $[16, 8, 4]$, we have more than $16$ code words of length $6$ which correspond to the equivalence classes of $\overline{X}$. For example, the following is an incident matrix of a $2$-$(16, 6, 2)$ block design from which we obtain $[16, 6, 6]$-linear code of $\mathbb{F}_2^{16}$:
\begin{equation}\label{des-inc}
{\small \left[\begin{array}{cccccccccccccccc}
 1 & 1 & 0 & 0 & 0 & 0 & 0 & 0 & 1 & 0 & 1 & 0 & 1 & 0 & 1 & 0 \\
 1 & 1 & 0 & 0 & 0 & 0 & 0 & 0 & 0 & 1 & 0 & 1 & 0 & 1 & 0 & 1 \\
 1 & 0 & 1 & 0 & 1 & 0 & 1 & 0 & 1 & 1 & 0 & 0 & 0 & 0 & 0 & 0 \\
 1 & 0 & 1 & 0 & 0 & 1 & 0 & 1 & 0 & 0 & 1 & 1 & 0 & 0 & 0 & 0 \\
 1 & 0 & 0 & 1 & 1 & 0 & 0 & 1 & 0 & 0 & 0 & 0 & 1 & 1 & 0 & 0 \\
 1 & 0 & 0 & 1 & 0 & 1 & 1 & 0 & 0 & 0 & 0 & 0 & 0 & 0 & 1 & 1 \\
 0 & 1 & 1 & 0 & 1 & 0 & 0 & 1 & 0 & 0 & 0 & 0 & 0 & 0 & 1 & 1 \\
 0 & 1 & 1 & 0 & 0 & 1 & 1 & 0 & 0 & 0 & 0 & 0 & 1 & 1 & 0 & 0 \\
 0 & 1 & 0 & 1 & 1 & 0 & 1 & 0 & 0 & 0 & 1 & 1 & 0 & 0 & 0 & 0 \\
 0 & 1 & 0 & 1 & 0 & 1 & 0 & 1 & 1 & 1 & 0 & 0 & 0 & 0 & 0 & 0 \\
 0 & 0 & 1 & 1 & 0 & 0 & 0 & 0 & 1 & 0 & 1 & 0 & 0 & 1 & 0 & 1 \\
 0 & 0 & 1 & 1 & 0 & 0 & 0 & 0 & 0 & 1 & 0 & 1 & 1 & 0 & 1 & 0 \\
 0 & 0 & 0 & 0 & 1 & 1 & 0 & 0 & 1 & 0 & 0 & 1 & 1 & 0 & 0 & 1 \\
 0 & 0 & 0 & 0 & 1 & 1 & 0 & 0 & 0 & 1 & 1 & 0 & 0 & 1 & 1 & 0 \\
 0 & 0 & 0 & 0 & 0 & 0 & 1 & 1 & 1 & 0 & 0 & 1 & 0 & 1 & 1 & 0 \\
 0 & 0 & 0 & 0 & 0 & 0 & 1 & 1 & 0 & 1 & 1 & 0 & 1 & 0 & 0 & 1
\end{array}\right].}
\end{equation}

Thus, we can determine the lattices $L$ for which $s_3 (L)$ is spherical $5$-design and $s_2 (L) = (A_1)^{16}$ uniquely up to isometry.

Now, let $\{ \varepsilon_i \}_{1 \leqslant i \leqslant 16}$ be another orthonormal basis of $\mathbb{R}^{16}$, and take an isometry which maps
\begin{equation}
\sqrt{2} e_{2 i - 1} \mapsto \varepsilon_{2 i - 1} + \varepsilon_{2 i}
\quad \text{and} \quad
\sqrt{2} e_{2 i} \mapsto \varepsilon_{2 i - 1} - \varepsilon_{2 i}
\quad \text{for every} \;
1 \leqslant i \leqslant 8.
\end{equation}
Then, this isometry leads the definition of $\Lambda_{16, 2, 1}$ in Section \ref{sec-def} from the above lattice which corresponds to the matrix (\ref{des-inc}). Actually, from the first and fourth row of the above matrix, we have the following correspondences:
\begin{align*}
\frac{e_1 + e_2 + e_9 + e_{11} + e_{13} + e_{15}}{\sqrt{2}} &\mapsto f_2 = \varepsilon_1 + \frac{\varepsilon_9 + \cdots + \varepsilon_{16}}{2},\\
\frac{e_1 + e_3 + e_5 + e_7 + e_9 + e_{10}}{\sqrt{2}} &\mapsto f_1 = \frac{\varepsilon_1 + \cdots + \varepsilon_8}{2} + \varepsilon_9.
\end{align*}\quad

\subsection{The cases of $s_2 (L) = (A_2)^8, (A_4)^4, (A_8)^2, A_{16}$}\quad

We can write $A_n = \{ \pm \varepsilon_i \mp \varepsilon_j \: ; \: 1 \leqslant i < j \leqslant n+1 \}$. We take $x_0 \in X$. If there exists $y \in A_n$ such that $(x_0, y) = 1$, then we may assume $y = \varepsilon_1 - \varepsilon_2$ without loss of generality. Furthermore, we write ${m_1}' := |\{ y \in A_n \: ; \: (x_0, y) = 1 \}|$. Now, we consider several distinct cases, namely $(i)$ $y = \pm (\varepsilon_1 - \varepsilon_2)$, $(ii)$ $y = \pm (\varepsilon_1 - \varepsilon_i), \pm (\varepsilon_2 - \varepsilon_i)$ for $3 \leqslant i \leqslant n+1$, and $(iii)$ $y = \pm (\varepsilon_i - \varepsilon_j)$ for $3 \leqslant i < j \leqslant n+1$.

For the second case $(ii)$, we have $(x_0, \varepsilon_1 - \varepsilon_i) - (x_0, \varepsilon_2 - \varepsilon_i) = (x_0, \varepsilon_1 - \varepsilon_2) = 1$. If we have $(x_0, \varepsilon_1 - \varepsilon_i) = -1$, then $(x_0, \varepsilon_2 - \varepsilon_i) = -2$, contradicting the relation (\ref{bd-inpro}). Thus, $(x_0, \varepsilon_1 - \varepsilon_i) = 0, 1$ and $(x_0, \varepsilon_2 - \varepsilon_i) = 0, -1$; just one of $(x_0, \varepsilon_1 - \varepsilon_i)$ and $(x_0, \varepsilon_2 - \varepsilon_i)$ is nonzero. Moreover, half of the vectors of this case are orthogonal to $x_0$.

For the third case $(iii)$, we have $(x_0, \varepsilon_i - \varepsilon_j) = - (x_0, \varepsilon_1 - \varepsilon_i) + (x_0, \varepsilon_1 - \varepsilon_j)$. $(x_0, \varepsilon_i - \varepsilon_j)$ is nonzero if and only if just one of $(x_0, \varepsilon_1 - \varepsilon_i)$ and $(x_0, \varepsilon_1 - \varepsilon_j)$ is equal to $1$. We put $D := |\{ 3 \leqslant i \leqslant n+1 \: ; \: (x_0, \varepsilon_1 - \varepsilon_i) = 1 \}|$, then $2 D (n - 1 - D)$ of vectors of this case are not orthogonal to $x_0$.

In conclusion, we have ${m_1}' = 1 + (n-1) + D (n-1 - D) = (n - D) (D + 1)$. Note that ${m_1}'$ is even if $n$ is even.

Recall that $|\{ y \in s_2 (L) \: ; \: (x_0, y) = 1 \}| = n_2$, and note that this number for $s_2 (L)$ must be a combination of the numbers ${m_1}'$ for each $A_n$. However, if $n = 2, 4, 8, 16$, then we have $n_2 = 9, 15, 27, 51$, where all of them are odd. These facts are contradictory. Thus, we can omit these cases, when $s_2 (L) = (A_n)^{16 / n}$ for $n > 1$.\\

\subsection{The cases of $s_2 (L) = (D_4)^4, (D_8)^2, D_{16}$}\quad

Following the procedure of the previous section, we write $D_n = \{ \pm \varepsilon_i \pm \varepsilon_j, \: \pm \varepsilon_i \mp \varepsilon_j \: ; \: 1 \leqslant i < j \leqslant n \}$ for an orthonormal basis of $\mathbb{R}^n$. We take $x_0 \in X$, then we assume $(x_0, \varepsilon_1 + \varepsilon_2) = 1$ if some vector of $D_n$ is not orthogonal to $x_0$. Furthermore, we write ${m_1}' := |\{ y \in D_n \: ; \: (x_0, y) = 1 \}|$. Now, we consider several distinct cases, namely $(i)$ $y = \pm (\varepsilon_1 + \varepsilon_2)$, $(ii)$ $y = \pm (\varepsilon_1 \pm \varepsilon_i), \pm (\varepsilon_2 \pm \varepsilon_i)$ for $3 \leqslant i \leqslant n$, $(iii)$ $y = \pm (\varepsilon_i \pm \varepsilon_j)$ for $3 \leqslant i < j \leqslant n$, and $(iv)$ $y = \pm (\varepsilon_1 - \varepsilon_2)$.

For the second case $(ii)$, we have $(x_0, \varepsilon_1 \pm \varepsilon_i) + (x_0, \varepsilon_2 \mp \varepsilon_i) = 1$. As in the section above, half of the vectors of this case are orthogonal to $x_0$.

For the fourth case $(iv)$, we put $D := (x_0, \varepsilon_1 + \varepsilon_i) + (x_0, \varepsilon_1 - \varepsilon_i) = (x_0, 2 \varepsilon_1)$ for any $3 \leqslant i \leqslant n$, which does not depend on $i$. If $D = 0$, then $(x_0, \varepsilon_2 + \varepsilon_i) + (x_0, \varepsilon_2 - \varepsilon_i) = 1$; thus we can take $\varepsilon_2$ instead of $\varepsilon_1$ without loss of generality. We may assume $D = 1, 2$.

If $D = 1$, then $(x_0, \varepsilon_1 - \varepsilon_2) = 0$. For the third case $(iii)$, we have $(x_0, \varepsilon_i \pm \varepsilon_j) = (x_0, \varepsilon_1 + \varepsilon_i) - (x_0, \varepsilon_1 \pm \varepsilon_j)$. Thus, just one of $(x_0, \varepsilon_i + \varepsilon_j)$ and $(x_0, \varepsilon_i - \varepsilon_j)$ is zero. Moreover, half of the vectors of the case $(iii)$ are orthogonal to $x_0$. In conclusion, we have ${m_1}' = 1 + 2 (n-2) + (n-2) (n-3) / 2 + 0 = n (n - 1) / 2$.

If $D = 2$, then $(x_0, \varepsilon_1 - \varepsilon_2) = 1$. For the third case $(iii)$, we have $(x_0, \varepsilon_i \pm \varepsilon_j) = (x_0, \varepsilon_1 + \varepsilon_i) - (x_0, \varepsilon_1 \pm \varepsilon_j) = 0$. Then, all of the vectors of the case $(iii)$ are orthogonal to $x_0$. In conclusion, we have ${m_1}' = 1 + 2 (n-2) + 0 + 1 = 2 (n - 1)$.

\begin{center}
\begin{tabular}{c}
\quad\\ \hline $D$\\ \hline ${m_1}'$\\ \hline
\end{tabular}
\quad
\begin{tabular}{cc}
\multicolumn{2}{c}{$(D_4)^4$}\\
\hline
$1$ & $2$\\
\hline
$6$ & $6$\\
\hline
\end{tabular}
\quad
\begin{tabular}{cc}
\multicolumn{2}{c}{$(D_8)^2$}\\
\hline
$1$ & $2$\\
\hline
$28$ & $14$\\
\hline
\end{tabular}
\quad
\begin{tabular}{cc}
\multicolumn{2}{c}{$D_{16}$}\\
\hline
$1$ & $2$\\
\hline
$120$ & $30$\\
\hline
\end{tabular}
\end{center}

Recall that $|\{ y \in s_2 (L) \: ; \: (x_0, y) = 1 \}| = n_2$, and note that this number must be the sum of the numbers ${m_1}'$. If $n = 4, 8, 16$, then we have $n_2 = 18, 42, 90$, respectively. Then, there remain some possibilities such that $n_2 = 6 + 6 + 6$ for the case of $s_2 (L) = (D_4)^4$ and $n_2 = 28 + 14$ for the case of $s_2 (L) = (D_8)^2$. On the other hand, $n_2 \ne {m_1}'$ for the case of $s_2 (L) = D_{16}$; thus we can omit this case.

Now, we write $x_0 = (a_1, \ldots, a_n, a_{n+1}, \ldots, a_{16})$, where each $a_i$ for $1 \leqslant i \leqslant n$ is the coordinate of $\varepsilon_i$ used for $D_n$. Then, we have $(x_0, \varepsilon_i \pm \varepsilon_j) = a_i \pm a_j$, thus we can calculate each $a_i$ from the values of $(x_0, \varepsilon_i \pm \varepsilon_j)$.

If $D = 1$, then we have $|a_i| = \frac{1}{2}$ for every $1 \leqslant i \leqslant n$. Furthermore, for $2 a_i \varepsilon_i + 2 a_j \varepsilon_j \in D_n$, $x_0 - (2 a_i \varepsilon_i + 2 a_j \varepsilon_j)$ is also a vector of norm $3$, where both signs of the $i$th and $j$th coordinates are different from those of $x_0$. Similarly, we have every element whose even signs are different from $x_0$. Here, we write
\begin{align*}
L_2 &:= \{(\pm 1, \ldots, \pm 1)/2 \in \mathbb{R}^n \, ; \, \text{the number of ``$-$'' is even}\},\\
L_3 &:= \{(\pm 1, \ldots, \pm 1)/2 \in \mathbb{R}^n \, ; \, \text{the number of ``$-$'' is odd}\},
\end{align*}
then the lattice $L$ includes either $L_2 \times \{ (a_{n+1}, \ldots, a_{16}) \}$ or $L_3 \times \{ (a_{n+1}, \ldots, a_{16}) \}$.

If $D = 2$, then we have $|a_i| = 1$ for some $1 \leqslant i \leqslant n$ and $a_j = 0$ for every $1 \leqslant j \leqslant n$ such that $i \ne j$. Furthermore, for $2 a_i \varepsilon_i$ and $a_i \varepsilon_i \pm \varepsilon_j \in D_n$, $x_0 - 2 a_i \varepsilon_i$ and $x_0 - (a_i \varepsilon_i \pm \varepsilon_j)$ is also a vector of norm $3$. Here, we write
\begin{equation*}
L_1 := \{ \pm \varepsilon_1, \ldots, \pm \varepsilon_n \},
\end{equation*}
then the lattice $L$ includes $L_1 \times \{ (a_{n+1}, \ldots, a_{16}) \}$.

On the other hand, if $x_0$ is orthogonal to $D_n$, then $a_i = 0$ for every $1 \leqslant i \leqslant n$. Here, we write
\begin{equation*}
L_0 := \{(0, \ldots, 0) \in \mathbb{R}^n \};
\end{equation*}
then we can say that the lattice $L$ includes $L_0 \times \{ (a_{n+1}, \ldots, a_{16}) \}$.

\subsubsection{The case of $s_2 (L) = (D_4)^4$}

We have $L_1 \simeq L_2 \simeq L_3 \simeq \mathbb{Z}^4$, where each vector is of norm $1$. Thus, $s_3 (L) \supset L_{i_1} \times L_{i_2} \times L_{i_3} \times L_{i_4}$, where $i_k = 0$ for just one $k$ and this set has $512$ vectors. Since $|s_3 (L)| = 2048$, $s_3 (L)$ includes four such sets. Then, we write
\begin{align*}
L_I &= (L_{i_{1,1}} \times L_{i_{1,2}} \times L_{i_{1,3}} \times L_{i_{1,4}})
 \cup (L_{i_{2,1}} \times L_{i_{2,2}} \times L_{i_{2,3}} \times L_{i_{2,4}})\\
 &\cup (L_{i_{3,1}} \times L_{i_{3,2}} \times L_{i_{3,3}} \times L_{i_{3,4}})
 \cup (L_{i_{4,1}} \times L_{i_{4,2}} \times L_{i_{4,3}} \times L_{i_{4,4}}),
\end{align*}
where $I = (i_{1,1}, \ldots, i_{4,4})$.

Firstly, for every $y_0 \in D_4$, we have $p_1 = 384$ and $|\{ x \in L_{i_{j,1}} \times L_{i_{j,2}} \times L_{i_{j,3}} \times L_{i_{j,4}} \: ; \: (x, y_0) = 1 \}| = 0$ or $128$. Thus, we have that just one of $i_{1,k}, i_{2,k}, i_{3,k}, i_{4,k}$ is zero for each $1 \leqslant k \leqslant 4$. We may assume that $i_{1,4} = i_{2,3} = i_{3,2} = i_{4,1} = 0$, without loss of generality. Secondly, we consider the inner product $(x_1, x_2)$ for $x_1 \in L_{i_k}$ and $x_2 \in L_{j_k}$. We have the fact that $(x_1, x_2) = 0, \pm 1$ if $i_k = j_k \ne 0$, that $(x_1, x_2) = \pm 1 / 2$ if $i_k \ne j_k$ and $i_k \cdot j_k \ne 0$, and that $(x_1, x_2) = 0$ if $i_k \cdot j_k = 0$. Thus, we have $i_{j_1,k} \ne i_{j_2,k}$ for $1 \leqslant k \leqslant 4$ if $j_1 \ne j_2$. Furthermore, $\{ i_{1,k}, i_{2,k}, i_{3,k}, i_{4,k} \} = \{ 0, 1, 2, 3 \}$ for $1 \leqslant k \leqslant 4$. Finally, there is an isometry which maps $L_i$ to $L_{i'}$ for every $1 \leqslant i \leqslant 4$ and some $1 \leqslant i' \leqslant 4$. Thus, for every pair $I = (i_{1,1}, \ldots, i_{4,4})$ which satisfies $i_{1,4} = i_{2,3} = i_{3,2} = i_{4,1} = 0$ and $\{ i_{1,k}, i_{2,k}, i_{3,k}, i_{4,k} \} = \{ 0, 1, 2, 3 \}$ for each $1 \leqslant k \leqslant 4$, $L_I$ is isometric to the following set:
\begin{align*}
&(L_2 \times L_2 \times L_1 \times L_0)
 \cup (L_3 \times L_3 \times L_0 \times L_1)\\
&\quad \cup (L_1 \times L_0 \times L_2 \times L_2)
 \cup (L_0 \times L_1 \times L_3 \times L_3)\\
&\subset (f_1 + (D_4)^4) \cup (f_1 - f_3 + (D_4)^4)\\
&\qquad \cup (f_2 + (D_4)^4) \cup (f_2 - f_3 + (D_4)^4).
\end{align*}

Thus, we can determine $s_3 (L)$ uniquely up to isometry, which generates $\Lambda_{16, 2, 2}$ in Section \ref{sec-def}.\\

\subsubsection{The case of $s_2 (L) = (D_8)^2$}

Every vector of $L_1$ is of norm $1$, and every vector of $L_2$ and $L_3$ is of norm $2$. Thus, $s_3 (L) \supset L_1 \times L_i$ or $L_j \times L_1$ for some $i, j = 2, 3$, where each set has $2048$ vectors. Similar to the previous section, we have
\begin{equation*}
L_{i, j} = (L_1 \times L_i) \cup (L_j \times L_1).
\end{equation*}
When we consider an isometry which changes the sign of one fixed coordinate of each vector, then this maps $L_2$ to $L_3$ and stabilizes $L_1$. Thus, for every $i, j = 2, 3$, $L_{i, j}$ is isometric to the following set:
\begin{equation*}
(L_2 \times L_1) \cup (L_1 \times L_2)
\; \subset \;
(f_1 + (D_8)^2) \cup (f_2 + (D_8)^2)
\end{equation*}

Thus we can determine $s_3 (L)$ uniquely up to isometry, which generates $\Lambda_{16, 2, 3}$ in Section \ref{sec-def}.\\

\subsection{The case of $s_2 (L) = (E_8)^2$}

Similar to the case of $s_2 (L) = (D_8)^2$, we write $E_8 = \{ \pm \varepsilon_i \pm \varepsilon_j, \: \pm \varepsilon_i \mp \varepsilon_j \: ; \: 1 \leqslant i < j \leqslant 8 \} \cup \{ (\pm \varepsilon_1 \pm \cdots \pm \varepsilon_8)/ 2 \: ; \: \text{the number of `$-$' is even} \}$ for an orthonormal basis of $\mathbb{R}^8$. We take $x_0 \in X$, then we assume that $(x_0, \varepsilon_1 + \varepsilon_2) = 1$ if such an element exists. Furthermore, we write ${m_1}' := |\{ y \in E_8 \: ; \: (x_0, y) = 1 \}|$. We consider several distinct cases; namely $(i)$ $y = \pm (\varepsilon_1 + \varepsilon_2)$, $(ii)$ $y = \pm (\varepsilon_1 \pm \varepsilon_i), \pm (\varepsilon_2 \pm \varepsilon_i)$ for $3 \leqslant i \leqslant 8$, $(iii)$ $y = \pm (\varepsilon_i \pm \varepsilon_j)$ for $3 \leqslant i < j \leqslant 8$, $(iv)$ $y = \pm (\varepsilon_1 - \varepsilon_2)$, $(v)$ $y = \pm (\varepsilon_1 + \varepsilon_2 \pm \varepsilon_3 \pm \cdots \pm \varepsilon_8) / 2$, and $(vi)$ $y = \pm (\varepsilon_1 - \varepsilon_2 \pm \varepsilon_3 \pm \cdots \pm \varepsilon_8) / 2$.

The cases $(i)$, $(ii)$, $(iii)$, and $(vi)$ are similar to the case of $s_2 (L) = (D_8)^2$. Note that $(x_0, \varepsilon_1 \pm \varepsilon_j) = 0, 1$ and $(x_0, \varepsilon_2 \pm \varepsilon_j) = 0, 1$, and we put $D := (x_0, \varepsilon_1 + \varepsilon_i) + (x_0, \varepsilon_1 - \varepsilon_i) \in \{ 1, 2 \}$.

For the fifth case $(v)$, we have $(x_0, \varepsilon_1 + \varepsilon_2) = (x_0, (\varepsilon_1 + \varepsilon_2 \pm \varepsilon_3 \pm \cdots \pm \varepsilon_8) / 2) + (x_0, (\varepsilon_1 + \varepsilon_2 \mp \varepsilon_3 \mp \cdots \mp \varepsilon_8) / 2)$. Thus, half of the vectors of this case are orthogonal to $x_0$.

If $D = 2$, we have $(x_0, \varepsilon_1 - \varepsilon_2) = 1$. The sixth case $(vi)$ is similar to the case $(v)$, and half of the vectors of this case are orthogonal to $x_0$. We have ${m_1}' = 1 + 12 + 0 + 1 + 16 + 16 = 46$.

If $D = 1$, just one of $(x_0, \varepsilon_i + \varepsilon_j)$ and $(x_0, \varepsilon_i - \varepsilon_j)$ is zero for the case $(iii)$. We may assume $(x_0, \varepsilon_3 + \varepsilon_4) = (x_0, \varepsilon_5 + \varepsilon_6) = (x_0, \varepsilon_7 + \varepsilon_8) = 1$ without loss of generality. Firstly, when the sign of $\varepsilon_3$ is equal to that of $\varepsilon_4$, then we have $(x_0, \varepsilon_3 + \varepsilon_4) = (x_0, (\pm \varepsilon_1 \mp \varepsilon_2 + \varepsilon_3 + \varepsilon_4 \pm \varepsilon_5 \pm \cdots \pm \varepsilon_8) / 2) - (x_0, (\pm \varepsilon_1 \mp \varepsilon_2 - \varepsilon_3 - \varepsilon_4 \pm \varepsilon_5 \pm \cdots \pm \varepsilon_8) / 2)$; thus half of such vectors are orthogonal to $x_0$. Secondly, when `the sign of $\varepsilon_3$ is not equal to that of $\varepsilon_4$' and `the sign of $\varepsilon_5$ is equal to that of $\varepsilon_6$', then $(x_0, \varepsilon_5 + \varepsilon_6) = (x_0, (\pm \varepsilon_1 \mp \varepsilon_2 \pm \varepsilon_3 \mp \varepsilon_4 + \varepsilon_5 + \varepsilon_6 \pm \varepsilon_7 \mp \varepsilon_8) / 2) - (x_0, (\pm \varepsilon_1 \mp \varepsilon_2 \pm \varepsilon_3 \mp \varepsilon_4 - \varepsilon_5 - \varepsilon_6 \pm \varepsilon_7 \mp \varepsilon_8) / 2)$, thus half of such vectors are orthogonal to $x_0$. Finally, when `the sign of $\varepsilon_3$ is not equal to that of $\varepsilon_4$' and `the sign of $\varepsilon_5$ is not equal to that of $\varepsilon_6$', then $(x_0, \varepsilon_7 + \varepsilon_8) = (x_0, (\pm \varepsilon_1 \mp \varepsilon_2 \pm \varepsilon_3 \mp \varepsilon_4 \pm \varepsilon_5 \mp \varepsilon_6 + \varepsilon_7 + \varepsilon_8) / 2) - (x_0, (\pm \varepsilon_1 \mp \varepsilon_2 \pm \varepsilon_3 \mp \varepsilon_4 \pm \varepsilon_5 \mp \varepsilon_6 - \varepsilon_7 - \varepsilon_8) / 2)$; thus half of such vectors are orthogonal to $x_0$. In conclusion, we have ${m_1}' = 1 + 12 + 15 + 0 + 16 + 16 = 60$.

Recall that $|\{ y \in s_2 (L) \: ; \: (x_0, y) = 1 \}| = n_2 = 90$, thus we cannot write $n_2$ as a sum of ${m_1}'$, so we can omit this case.\\

\begin{remark}
In the argument of Section \ref{sec-16-2-1}, from the lattice $\Lambda_{16, 2, 1}$, we obtain a $2$-$(16, 6, 2)$ block design which generates $[16, 6, 6]$-linear code of $\mathbb{F}_2^{16}$. Conversely, it determines the lattice $\Lambda_{16, 2, 1}$ uniquely.

On the other hand, there are two more equivalence classes of $2$-$(16, 6, 2)$ block designs. Note that we obtain the lattices $\Lambda_{16, 2, 2}$ and $\Lambda_{16, 2, 3}$ from the two equivalence classes which generate $[16, 7, 4]$ and $[16, 8, 4]$-linear codes of $\mathbb{F}_2^{16}$, respectively.
\end{remark}\quad

\section*{Acknowledgment}
Many thanks to Professor Eiichi Bannai for suggesting these problems. We thank the anonymous referee for many useful suggestions and comments.\\

\end{document}